\documentclass[11pt]{article}
\usepackage{latexsym,amsmath,amsfonts,amssymb,amstext,amsthm}

\newtheorem{teo}{Theorem}[section]%[section]
\newtheorem{coro}[teo]{Corollary}
\newtheorem{pro}[teo]{Proposition}
\newtheorem{lem}[teo]{Lemma}
\theoremstyle{definition}
\newtheorem{defi}{Definition}[section]
\newtheorem{rem}[defi]{Remark}
\newcommand{\prue}{\noindent \textbf{Proof:~~}}
\newcommand{\R}{\mathbb{R}} %% Conjunto de los N\'{u}meros Reales
\newcommand{\N}{\mathbb{N}} %Conjuntos Naturales
\newcommand{\Be}{\mathcal{B}}

\newcommand{\fin}{ \hfill $\square$}

\newcommand{\tlim}{\displaystyle\lim}

\newcommand{\defrac}{\displaystyle\frac}

\DeclareMathOperator{\Int}{Int}
\title{Common fixed points for a pair of commuting mappings in complete cone metric spaces}

\author{Jos\'e R.~Morales and Edixon M.~Rojas\\ \begin{scriptsize}moralesj@ula.ve\hspace{.7cm} edixonr@ula.ve\end{scriptsize} \\  \begin{scriptsize}Departamento de Matem\'aticas, Universidad de Los Andes\end{scriptsize}\\\begin{scriptsize}5101, M\'erida--Venezuela\end{scriptsize}} 

\date{}
\begin{document}
\maketitle

\begin{abstract}

This paper is devoted to prove the S. L. Singh's common fixed point Theorem for commuting mappings in cone metric spaces. In this framework, we introduce the notions of Generalized Kannan Contraction, Generalized Zamfirescu Contraction and Generalized Weak Contraction for a pair of mappings, proving afterward their respective fixed point results.
\end{abstract}

\section{Introduction and Preliminaries.}

In 1977, S.L. Singh \cite{Si} proved the following result.
\begin{teo}\label{Singh}

Let $(M,d)$ be a complete metric space. Let $S$ and $T$ be mappings from $M$ into itself such that,

\begin{itemize}

\item[(a)] $T$ is continuous;

\item[(b)] $S(M)\subset T(M);$

\item[(c)] $S$ and $T$ commute;

\item[(d)] The following inequality holds
 \begin{align*}
d(Sx,Sy)\leq& a d(Tx,Ty)+b[d(Sx,Ty)+d(Sx,Sy)]\\
&+c[d(Sx,Ty)+d(Sy,Tx)]
\end{align*}
 for all $x,y\in M,$  where $a,b,c$ are nonnegative real numbers such that $0<a+2b+2c<1.$ Then $S$ and $T$ have a unique common fixed point in $M.$
\end{itemize}
\end{teo}
As a consequence of this theorem, the results of Jungck (\cite{Ju}) can be obtained by considering the particular case $a=b=0$. On the other hand, when in Theorem \ref{Singh} is considered $T=Id$, the identity map, then are obtained the results by the Hardy \& Rogers  in \cite{HaRo73} .

The main goal of this paper is to present the Singh's fixed point Theorem \ref{Singh} in the setting of complete cone metric spaces.
Furthermore, we will introduce some contractive conditions for a pair of mappings on these spaces which generalize some well-known notions given in complete metric spaces.

First, we must recall that the cone metric spaces were introduced in 2007 by Huang and Zhang in \cite{GX}. They also obtained several fixed point theorems for contractive single valued maps in such spaces. Since then, a lot of works in this subject were already published, including various coincidence and common fixed point theorems for a pair of weakly compatible mapping (\cite{AJ}) as well as generalized contraction and Zamfirescu pair (\cite{BaAlVa10}).

\begin{defi}[\cite{GX}]
Let $(E,\|\cdot\|)$ be a real Banach space. A subset $P\subset E$ is called a cone if and only if:

\begin{description}

\item[(P1)] $P$ is closed, non empty and $P\neq \{0\}.$

\item[(P2)] $a,b\in \R,\quad a,b\geq 0,\quad x,y\in P$ implies $ax+by\in P.$

\item[(P3)] $x\in P$ and $-x\in P\Rightarrow x=0,$ that is, $P\cap (-P)=\{0\}.$
\end{description}
\end{defi}

Given a cone $P\subset E,$ we define a partial ordering $\leq$ with respect to $P$ by
\begin{equation*}
x\leq y\quad \mbox{if and only if}\quad y-x\in P.
\end{equation*}
 We write $x<y$ to indicate that $x\leq y$ but $x\neq y,$ while $x\ll y$ will stand for $y-x\in \Int\, P,$ where $\Int\, P$ denote the interior of $P$.
\begin{defi}[\cite{GX}]
Let $(E,\|\cdot\|)$ be a Banach space and $P\subset E$ a cone. The cone $P$ is a called normal if there is a number $K>0$ such that for all $x,y\in M$

$$0\leq x\leq y\quad \mbox{implies}\quad \|x\|\leq K\|y\|.$$ The least positive number satisfying the above is called the normal constant of $P.$
\end{defi}

In the following, we always suppose that $(E,\|\cdot\|)$ is a real Banach space, $P$ is a cone in $E$ with $\Int\, P\neq \emptyset$ and $\leq$ is partial ordering with respect to $P.$

\begin{defi}[\cite{GX}]
Let $M$ be a nonempty set. Suppose that the mapping $d: M\times M\longrightarrow E$ satisfies:
\begin{description}
\item[(CM1)] $0<d(x,y)$ for all $x,y\in M$ and $d(x,y)=0$ if and only if $x=y;$
\item[(CM2)] $d(x,y)=d(y,x)$ for all $x,y\in M;$
\item[(CM3)] $d(x,y)\leq d(x,z)+d(z,y)$ for all $x,y,z\in M.$
\end{description}
Then, $d$ is called a cone metric on $M,$ and the pair $(M,d)$ is called a cone metric space. It will be denoted by CMS.

\end{defi}
Note that the notion of cone metric space (CMS) is more general that the concept of metric space.

\begin{defi}[\cite{GX}]
Let $(M,d)$ be a CMS. Let $(x_n)$ be a sequence in $M$ and $x\in M.$
\begin{itemize}
\item[(i)] $(x_n)$ is said convergent to $x$ whenever for every $c\in E,$ with $0\ll c$ there is a
           positive integer $n_0$ such that $d(x_n,x)\ll c$ for all $n\geq n_0.$ We denote this by $\tlim_{n\rightarrow \infty} x_n=x$ or $x_n\rightarrow x$ as $n\rightarrow \infty.$
\item[(ii)] $(x_n)$ is said to be a Cauchy sequence in $M$ whenever for every $c\in E$ with $0\ll c$
            there is a positive integer $n_0$ such that $d(x_n,x_m)\ll c$ for all $n,m\geq n_0.$

\item[(iii)] $(M,d)$ is called a complete CMS if every Cauchy sequence is convergent in $M.$

\item[(iv)] A set $A\subseteq M$ is said to be closed if for any sequence $(x_n)\subset A$ convergent
            to $x,$ we have that $x\in A.$

\item[(v)] A set $A\subseteq M$ is called sequentially compact if for any sequence $(x_n)\subset A,$
           there exists a subsequence $(x_{n_k})$ of $(x_n)$ which is convergent to an element of $A$.
\end{itemize}
\end{defi}

\begin{lem}[\cite{GX}]

Let $(M,d)$ be a CMS, $P\subset E$ a normal cone with normal constant $K.$ Let $(x_n)$  be a sequence in $M$ and $x,y\in M.$
\begin{itemize}
\item[(i)] $(x_n)$ converges to $x$ if and only if $\tlim_{n\rightarrow \infty} d(x_n,x)=0$.

\item[(ii)] If $(x_n)$ converges to $x$, and $(x_n)$ converges to $y$, then $x=y.$

\item[(iii)] If $(x_n)$ converges to $x,$ then $(x_n)$ is a Cauchy sequence.

\item[(iv)] $(x_n)$ is a Cauchy sequence if and only if $\tlim_{n,m\rightarrow \infty} d(x_n,x_m)=0$.
\end{itemize}
\end{lem}

\begin{lem}[\cite{LS}]\label{Cauchy cond}
Let $(M,d)$ be a CMS and let $P\subset E$ be a normal cone with normal constant $K.$ If there exists a sequence $(x_n)$ in $M$ and $a$ real number $a\in (0,1)$ such that for every $n\in \N,$
\begin{equation*}
d(x_{n+1},x_n)\leq a d(x_n,x_{n-1}),
\end{equation*}
then $(x_n)$ is a Cauchy sequence.
\end{lem}

\begin{defi}[\cite{TA}]
Let $(M,d)$ be a CMS and $A\subset M.$
\begin{itemize}
\item[(i)] A point $b\in A$ is called an interior point of $A$ whenever there exists a point $c,\,\,
           0\ll c$ such that

           $$B(b,c)\subseteq A$$ where $B(b,c)=\{y\in M\,:\, d(y,b)\ll c\}.$

\item[(ii)] A subset $A\subset M$ is called open if each element of $A$ is an interior point of $A.$
\end{itemize}

The family $\Be=\{B(b,c)\,/\, b\in M,\,\, 0\ll c\}$ is a sub-basis for a topology on $M.$ We denote this cone topology by $\tau_c.$
\end{defi}

The topology $\tau_c$ is Hausdorff and first countable, (\cite{TA}). Hence, we conclude that any CMS $(M,d)$ is Hausdorff and the limits are unique.

\begin{defi}[\cite{TA}]
Let $(M,d)$ be a CMS. A map $T: M\longrightarrow M$ is called continuous at $x\in M,$ if for each $V\in \tau_c$ containing $Tx,$ there exists $U\in \tau_c$ containing $x$ such that
\begin{equation*}
T(U)\subset V.
\end{equation*}
If $T$ is continuous at each $x\in M$, then it is called continuous.
\end{defi}

\begin{defi}[\cite{TA}]
Let $(M,d)$ be a CMS. A map $T: M\longrightarrow M$ is called sequentially continuous if $(x_n)\subset M,\,\, x_n\rightarrow x$ implies $Tx_n\longrightarrow Tx.$
\end{defi}

\begin{pro}[\cite{TA}]\label{seque cont}
Let $(M,d)$ be a CMS and $T: M\longrightarrow M$ be any map. Then, $T$ is continuous if and only if $T$ is sequentially continuous.
\end{pro}

\section{On the Singh's common fixed point Theorem for commuting mappings in cone metric spaces}

In this section we will prove the Singh's common fixed point Theorem for commuting mappings in the framework of cone metric spaces. Afterwards, we are going to give some consequences of this result.

\begin{defi}[\cite{Ve}]

Let $(M,d)$ be CMS, and $P\subset E$ a normal cone with normal constant $K.$ Let $S,T: M\longrightarrow M$ be mappings such that $S(M)\subset T(M)$ and for every $x_0\in M$ we define the sequence $(x_n)$ by $T(x_n)=S(x_{n-1}),$ $n=1,2,\ldots,$ we say that $S(x_n)$ is a $(S,T)-$sequence with initial point $x_0.$
\end{defi}
%
%\begin{defi}[e.g.,~\cite{BaVe09}]
%Let $S$ and $T$ be self-maps of a set $M.$ If $w=S(z)=T(z)$ for some $z\in M,$ then $z$ is called a coincidence point of $S$ and $T,$ and $w$ is called a point of coincidence of $S$ and $T$. The mappings $S$ and $T$ are weakly compatible if, for every $z\in M$ holds:
%If $Sz=Tz$, then $STz=TTz.$
%\end{defi}

%\begin{pro}[\cite{AJ}]
%Let $S$ and $T$ be weakly compatible self-maps a set $M.$ If $S$ and $T$ have a unique point of coincidence $w=Sz=Tz$ then $w$ is the unique common fixed point of $S$ and $T.$
%\end{pro}

\begin{teo}\label{teo Singh cone}
Let $(M,d)$ be a complete CMS, and $P\subset E$ a normal cone with normal constant $K.$ Let $S$ and $T$ be self-maps of $M$ such that,
\begin{itemize}

\item[(a)] $T$ is continuous.

\item[(b)] $S(M)\subset T(M).$

\item[(c)] $(S,T)$ is a commuting pair.

\item[(d)] The following inequality holds
\begin{align*}
d(Sx,Sy)\leq& a d(Tx,Ty)+b[d(Sx,Tx)+d(Sy,Ty)]\\
&+c[d(Sx,Ty)+d(Sy,Ty)]\qquad\qquad\qquad\qquad\qquad\quad(S)
\end{align*}
for all $x,y\in M,$ where where $a,b,c$ are nonnegative real numbers such that $0<a+2b+2c<1.$
\end{itemize}
 Then $S$ and $T$ have a unique common fixed point.
\end{teo}
\prue Suppose that $x_0\in M$ is an arbitrary point. We will prove that the $(S,T)-$sequence $(S(x_n))$ with initial point $x_0$ is a Cauchy sequence in $M.$
In fact, notice that
\begin{align*}
d(Sx_{n+1},Sx_n)\leq& a d(Tx_{n+1},Tx_n)+b[d(Sx_{n+1},Tx_{n+1})+d(Sx_n,Tx_n)]\\
&+c[d(Sx_{n+1},Tx_n)+d(Sx_n,Tx_{n})].
\end{align*}
Thus,
\begin{equation*}
d(Sx_{n+1},Sx_n)\leq \defrac{a+b+c}{1-b-c}d(Sx_n,Sx_{n-1})
\end{equation*}
 which is the same that
\begin{equation}\label{eq2.1}
d(Sx_{n+1},Sx_n)\leq \alpha\, d(Sx_n,Sx_{n-1})
\end{equation}
with
\begin{equation*}
\alpha=\defrac{a+b+c}{1-b-c}<1.
\end{equation*}
In this way, from inequality \eqref{eq2.1} and Lemma \ref{Cauchy cond}, we have that $(S(x_n))$ is a Cauchy sequence in $M.$

On the other hand, repeating the procedure above we can conclude that for all $n\in \N$
\begin{equation*}
d(Sx_{n+1},Sx_n)\leq \alpha^n\, d(Sx_1,Sx_0),
\end{equation*}
since $P\subset E$ is a normal cone with normal constant $K$, then we have
\begin{equation*}
\|d(Sx_{n+1},Sx_n)\|\leq \alpha^n\, K\|d(Sx_1,Sx_0)\|.
 \end{equation*}
 Taking limits in inequality above we conclude that 
 \begin{equation*}
 \tlim_{n\rightarrow \infty} d(Sx_{n+1},Sx_n)=0.
 \end{equation*}
Since $M$ is a complete CMS, then there exists $z_0\in M$ such that
\begin{equation}\label{eq2.2}
\tlim_{n\rightarrow \infty} Sx_n=\tlim_{n\rightarrow \infty} Tx_{n+1}=z_0.
\end{equation} Since $T$ is continuous (and therefore sequentially continuous by Proposition \ref{seque cont}), also due to the fact that $S$ and $T$ commute, we have
\begin{equation*}
Tz_0=T\left(\tlim_{n\rightarrow \infty} Tx_n\right)=\tlim_{n\rightarrow \infty} T^2x_n
\end{equation*}
 as well as
\begin{equation*}
Tz_0=T\left(\tlim_{n\rightarrow \infty} Sx_n\right)=\tlim_{n\rightarrow \infty}TSx_n=\tlim_{n\rightarrow \infty} STx_n.
\end{equation*}
Now,
\begin{align*}
d(STx_n,Sz_0)&\leq a d(T^2x_n,Tz_0)+b[d(STx_n,T^2x_n)+d(Sz_0,Tz_0)]\\
&+c[d(STx_n,Tz_0)+d(Sz_0,Tz_0)].
\end{align*}
 Again, since $P$ is a normal cone with normal constant $K$ we have
\begin{align*}
\|d(STx_n,Sz_0)\|\leq& K[a \|d(T^2x_n,Tz_0)\|+b\|d(STx_n,T^2x_n)\|\\
&+c\|d(STx_n,Tz_0)\|+(b+c)\|d(Sz_0,Tz_0)\|]
\end{align*}
 taking the limit as $n\rightarrow \infty$ we obtain
\begin{align*}
\|d(Tz_0,Sz_0)\|\leq& K[a \|d(Tz_0,Tz_0)\|+b\|d(Tz_0,Tz_0)\|+c\|d(Tz_0,Tz_0)\|\\
&+(b+c)\|d(Sz_0,Tz_0)\|]
\end{align*}
or, rewriting the inequality above,
\begin{equation*}
\|d(Tz_0,Sz_0)\|\leq K(b+c)\|d(Sz_0,Tz_0)\|.
\end{equation*}
Hence, since $0\leq b+c<1,$ we have then $d(Tz_0,Sz_0)=0,$ that is $Tz_0=Sz_0.$

Now,
\begin{align*}
d(Sx_n,Sz_0)\leq& a d(Tx_n,Tz_0)+b[d(Sx_n,Tx_n)+d(Sz_0,Tz_0)]\\
&+c[d(Sx_n,Tz_0)+d(Sz_0,Tz_0)]\\
=&a d(Tx_n,Tz_0)+bd(Sx_n,Tx_n)+c d(Sx_n,Tz_0)\\
&+(b+c)d(Sz_0,Tz_0),
\end{align*}
since $P\subset E$ is a normal cone with normal constant $K,$ then we have
\begin{align*}
\|d(Sx_n,Sz_0)\|&\leq K[a \|d(Tx_n,Tz_0)\|+b\|d(Sx_n,Tx_n)\|+c\|d(Sx_n,Tz_0)\| \\
&+(b+c)\|d(Sz_0,Tz_0)\|].
\end{align*}
 Again, taking the limit as $n\rightarrow \infty$ we obtain,
\begin{align*}
\|d(z_0,Sz_0)\|\leq& K[a \|d(z_0,Tz_0)\|+b\|d(z_0,z_0)\|+c\|d(z_0,Tz_0)\|\\
&+(b+c)\|d(Sz_0,Tz_0)\|]\\
=& K(a+c)\|d(z_0,Tz_0)\|.
 \end{align*}
As above, we conclude that $d(z_0,Sz_0)=0$, which implies that $z_0=Sz_0$ and thus we have proved that
\begin{equation*}
Sz_0=Tz_0=z_0.
\end{equation*}

The uniqueness of the common fixed point $z_0$ follows from inequality $(S)$. In fact, let us suppose that $y_0=Sy_0=Ty_0$. Then,
\begin{align*}
d(y_0,z_0)=&d(Sy_0,Sz_0) \leq a d(Ty_0,Tz_0)+b[d(Sy_0,Ty_0)+d(Sz_0,Tz_0)]\\
&+c[d(Sy_0,Tz_0)+d(Sz_0,Tz_0)]\\
=&(a+c)d(Sy_0,Sz_0)=(a+c)d(y_0,z_0).
\end{align*}
As before, the conclusion follows from the fact that $0\leq a+c<1$. Thus the theorem is proved.\fin

\subsection{Some consequences of Theorem \ref{teo Singh cone}}

In this part we are going to mention some results, which now can be obtained as a consequence of Theorem  \ref{teo Singh cone}.
First, notice that if in Theorem \ref{teo Singh cone} we take $E=\R_{+}$ and $P=[0,+\infty)$ we obtain Theorem \ref{Singh} for a pair $(S,T)$ of mappings.
Now if we take $b=c=0$ in inequality $(S)$, we obtain the following.

\begin{coro}\label{Jongck}
Let $(M,d)$ be a complete CMS, and $P\subset E$ a normal cone with normal constant $K.$ Let $S$ and $T$ be self-maps of $M$ such that,
\begin{itemize}
\item[(a)] $T$ is continuous.

\item[(b)] $S(M)\subset T(M).$

\item[(c)] $(S,T)$ is a commuting pair.

\item[(d)] The inequality
\begin{equation*}
d(Sx,Sy)\leq a d(Tx,Ty)\qquad\qquad\qquad \qquad\qquad\qquad\qquad (J)
\end{equation*}
holds for all $x,y\in M$, where $0\leq a<1$.
\end{itemize}
 Then $S$ and $T$ have a unique common fixed point.
\end{coro}
In this case, if we take in the Corollary \ref{Jongck}, $E=\R_{+}$ and $P=[0,+\infty)$, then we obtain the result given in 1976 by G. Jungck in \cite{Ju}.
On the other hand, if we consider $a=c=0$ in $(S)$, then we obtain the next result.
\begin{coro}\label{Kannan}
Let $(M,d)$ be a complete CMS and $P\subset E$ a normal cone with normal constant $K.$ Let $S$ and $T$ be self maps of $M$ such that,
\begin{itemize}

\item[(a)] $T$ is continuous.

\item[(b)] $S(M)\subset T(M).$

\item[(c)] $(S,T)$ is a commuting pair.

\item[(d)] The inequality
\begin{equation*}
d(Sx,Sy)\leq b[d(Sx,Tx)+d(Sy,Ty)]\qquad \qquad \qquad \qquad(GKC)
\end{equation*}
 is satisfies for all $x,y\in M$ and $0\leq b<1/2$.
\end{itemize}
Then $S$ and $T$ have a unique common fixed point.
\end{coro}

 The mappings satisfying inequality (GKC) are called \emph{generalized Kannan contractions}.
 If in the Corollary \ref{Kannan} we take $E=\R_{+}$, $P=[0,+\infty)$ and $T=Id$ (the identity mapping) we obtain the Kannan's result \cite{Ka}.
Finally, if we consider $a=b=0$ in inequality $(S)$, then we obtain the following result.

\begin{coro}\label{chatterjea}
Let $(M,d)$ be a complete CMS, and $P\subset E$ a normal cone with normal constant $K.$ Let $S$ and $T$ be self-maps of $M$ such that,

\begin{itemize}
\item[(a)] $T$ is continuous.

\item[(b)] $S(M)\subset T(M).$

\item[(c)] $(S,T)$ is a commuting pair.

\item[(d)] the inequality
\begin{equation*}
d(Sx,Sy)\leq c[d(Sx,Ty)+d(Sy,Tx)]\qquad \qquad\qquad \qquad (GCC)
\end{equation*}
holds for all $x,y\in M$ and $0\leq c<1/2.$
\end{itemize}
Then $S$ and $T$ have a unique common fixed point.
\end{coro}

Notice that if in the Corollary \ref{chatterjea}, we take $E=\R_{+}$, $P=[0, +\infty)$ and $T=Id$, then we obtain the Chatterjea's results \cite{Cha}.
 The mappings satisfying condition (GCC) are called \emph{generalized Chatterjea contractions}.

\section{Common fixed points for generalized Zamfrescu and weak contraction operators on CMS}

Using the ideas of T. Zamfrescu \cite{Za} (see also, \cite{OI}) we introduce the notion of Generalized Zamfirescu operators (GZ0) in the framework of complete cone metric spaces.

\begin{defi}\label{def zamfirescu}
Let $(M,d)$ be a CMS and let $S,T: M\longrightarrow M$ be two mappings. The pair $(S,T)$ is called the generalized Zamfirescu operators, (GZ0), if there are $0\leq a<1,\,\, 0\leq b<1/2$ and $0\leq c<1/2$ such that for all $x,y\in M,$ at least one of the next conditions are true:

\begin{description}

\item[(GZ01)] $d(Sx,Sy)\leq a d(Tx,Ty)$

\item[(GZ02)] $d(Sx,Sy)\leq b[d(Sx,Tx)+d(Sy,Ty)]$

\item[(GZ03)] $d(Sx,Sy)\leq c[d(Sx,Ty)+d(Sy,Tx)]$
\end{description}
\end{defi}

\begin{rem}$\mbox{}$
\begin{enumerate}
\item If in Definition \ref{def zamfirescu} we take $M$ a Banach space, $E=\R_{+}$ and $P=[0, +\infty)$, then we
      obtain the definition given by M. O. Olantinwo and C. O. Imoru \cite{OI}.

\item If in Definition \ref{def zamfirescu} we take $E=\R_{+},\,\, P=[0, +\infty)$ and $T=Id,$ identity mapping,
      then we get the Zamfirescu's definition \cite{Za}.
\end{enumerate}
\end{rem}

From Definition \ref{def zamfirescu} is immediate the following.

\begin{pro}\label{equiv zamfirescu}
Let $(M,d)$ be a CMS and $S,T: M\longrightarrow M$ a pair of \textup{(GZ0)}, then we have
\begin{itemize}
\item[(a)] $d(Sx,Sy)\leq \delta d(Tx,Ty)+2\delta d(Sx,Tx)$.

\item[(b)] $d(Sx,Sy)\leq \delta d(Tx,Ty)+2\delta d(Sy,Tx)$.
\end{itemize}
for all $x,y\in M$ and where
\begin{equation*}
\delta=\max\left\{a,\, \defrac{b}{1-b},\, \defrac{c}{1-c}\right\},\quad 0\leq\delta<1.
\end{equation*}
\end{pro}

The following result generalize the well--known theorem given by T. Zamfirescu in \cite{Za}.

\begin{teo}\label{GZ0 teo}

Let $(M,d)$ be a complete CMS, and $P\subset E$ a normal cone with normal constant $K.$ Suppose that $S,T: M\longrightarrow M$ are \textup{(GZ0)} such that,

\begin{itemize}
\item[(a)] $T$ is continuous.

\item[(b)] $S(M)\subset T(M).$

\item[(c)] $(S,T)$ is a commuting pair.
\end{itemize}
 Then, $S$ and $T$ have a unique common fixed point.
\end{teo}
\prue
Since the pair $(S,T)$  is a (GZ0), then by Proposition \ref{equiv zamfirescu} (b) we have that
\begin{equation*}
d(Sx,Sy)\leq \delta d(Tx,Ty)+2\delta d(Sy,Tx),\qquad\forall x,y\in M,
\end{equation*}
 where
 \begin{equation*}
1>\delta=\max\left\{a,\, \defrac{b}{1-b},\, \defrac{c}{1-c}\right\}.
\end{equation*}
Suppose that $x_0\in M$ is an arbitrary point. We are going to prove that the $(S,T)$-sequence $S(x_n)$ with initial point $x_0$ is a Cauchy sequence in $M.$

Notice that
\begin{align*}
d(Sx_{n+1},Sx_n)\leq& \delta d(Tx_{n+1},T_nx)+2\delta d(Sx_{n},Tx_{n+1})\\
=&\delta d(Sx_n,Sx_{n-1})+2\delta d(Sx_{n},Sx_n)\\
=&\delta d(Sx_n,Sx_{n-1}).
\end{align*}
Therefore,  from Lemma \ref{Cauchy cond} we conclude that $S(x_n)$ is a Cauchy sequence.
Repeating the procedure above, we get
\begin{equation*}
d(Sx_{n+1},Sx_n)\leq\delta^nd(Sx_1,Sx_0).
\end{equation*}
Thus, taking norm we obtain
\begin{equation*}
\|d(Sx_{n+1},Sx_n)\|\leq\delta^nK\|d(Sx_1,Sx_0)\|
\end{equation*}
therefore, as the proof of Theorem \ref{teo Singh cone}, taking the limit as $n\to\infty$ and using the fact that $M$ is a complete CMS, we guarantee the existence of a $z_0\in M$ such that
\begin{equation}\label{Gen Zam eq1}
\lim_{n\to \infty}Sx_n=\lim_{n\to \infty}Tx_{n+1}=z_0.
\end{equation}
By the continuity of $T$ and the commutating between $T$ and $S$, we have
\begin{align}
Tz_0=&\lim_{n\to\infty}T^2x_n \label{Gen Zam eq2}\\
Tz_0=&\lim_{n\to\infty}STx_n. \label{Gen Zam eq3}
\end{align}
On the other hand, using Proposition \ref{equiv zamfirescu} (a), we get
\begin{equation*}
d(STx_n,Sz_0)\leq\delta d(T^2x_n,Tz_0)+2\delta d(STx_n,T^2x_n).
\end{equation*}
Taking norm and the limit as $n\to\infty$ to the above inequality we have
\begin{equation*}
\lim_{n\to\infty}\|d(STx_n,Sz_0)\|\leq K\lim_{n\to\infty}[\delta\|d(T^2x_n,Tz_0)\|+2\delta\|d(STx_n,T^2x_n)\|].
\end{equation*}
From \eqref{Gen Zam eq2} and \eqref{Gen Zam eq3} we conclude that
\begin{equation*}
d(Tz_0,Sz_0)=0
\end{equation*}
that is, $Tz_0=Sz_0$.

Now, using Proposition \ref{equiv zamfirescu} (a), the following inequality holds,
\begin{equation*}
d(Sx_n,Sz_0)\leq\delta d(Tx_n,Tz_0)+2\delta d(Sx_n,Tx_n).
\end{equation*}
Repeating the argument above we obtain
\begin{equation*}
\lim_{n\to\infty}\|d(Sx_n,Sz_0)\|\leq K\lim_{n\to\infty}[\delta\|d(Tx_n,Tz_0)\|+2\delta\|d(Sx_n,Tx_n)\|],
\end{equation*}
using \eqref{Gen Zam eq1} and the fact that $Tz_0=Sz_0$ we conclude that
\begin{equation*}
\|d(z_0,Sz_0)\|\leq K\delta\|d(z_0,Sz_0)\|.
\end{equation*}
Since $K\delta<1$, then $z_0=Sz_0$.

Finally, we are going to prove the uniqueness of the fixed point. Let us suppose that $y_0\in M$ is such that $y_0=Sy_0=Ty_0$. Then
\begin{align*}
d(z_0,y_0)=d(Sz_0,Sy_0)\leq&\delta d(Tz_0,Ty_0)+2\delta d(Sz_0,Tz_0)\\
=&\delta d(z_0,y_0).
\end{align*}
From the fact that $\delta<1$, we conclude that $z_0=y_0$. Thus the theorem is proved.
 \fin

Now we introduce an equivalent notion of (GZ0) in CMS as follows.

\begin{defi}\label{def zamfirescu equi}
Let $(M,d)$ be a CMS and $S,T: M\longrightarrow M$ be mappings for which there exists $0\leq h<1$ such that for all $x,y\in M$
\begin{align}\label{def zamfirescu equi equa}
&d(Sx,Sy)\leq&\nonumber\\
 &h\max\left\{d(Tx,Ty),\, \frac{d(Sx,Tx)+d(Sy,Ty)}{2},\, \frac{d(Sx,Ty)+d(Sy,Tx)}{2}\right\}.
 \end{align}
 \end{defi}
It is not difficult to see that Definitions \ref{def zamfirescu} and \ref{def zamfirescu equi} are equivalent. Therefore, our results remain valid for mappings satisfying \eqref{def zamfirescu equi equa} as well.

In 2003, V. Berinde, (\cite{Be} and \cite{Ber}), introduced a new class of contraction mappings on metric spaces, which are called \emph{weak contraction}. Now we will extend these kind of contractive condition to a pair of mappings in the setting of CMS.

\begin{defi}
Let $(M,d)$ be a CMS and $S,T: M\longrightarrow M$ two mappings. The pair $(S,T)$ is called a generalized weak contraction, (GWC) if there exist constants $0<\delta<1$ and $L\geq 0$ such that
\begin{equation*}
d(Sx,Sy)\leq \delta d(Tx,Ty)+L d(Sx,Tx)
\end{equation*}
 for all $x,y\in M$.
\end{defi}

The next proposition gives examples of (GWC). This can be proved in a similar form as Proposition 3.3 of \cite{MoRo}.

\begin{pro}\label{GWC prop}

Let $(M,d)$ be a CMS and $S,T: M\longrightarrow M$ two mappings. Then,
\begin{itemize}
\item[(a)] If the pair $(S,T)$ satisfies the condition $(J)$, then $(S,T)$ is a \textup{(GWC)}.

\item[(b)] If the pair $(S,T)$ is a \textup{(GKC)}, then $(S,T)$ is a \textup{(GWC)}.

\item[(c)] If the pair $(S,T)$ is a \textup{(GCC)}, then $(S,T)$ is a \textup{(GWC)}.

\item[(d)] If the pair $(S,T)$ is a \textup{(GZ0)}, then $(S,T)$ is a \textup{(GWC)}.
\end{itemize}
\end{pro}
Next result can be proved by following the procedure used in the proof of Theorem \ref{GZ0 teo}.
\begin{teo}
Let $(M,d)$ be a complete CMS and $S,T: M\longrightarrow M$ two mappings. If the pair $(S,T)$ is a \textup{(GWC)} such that,
\begin{itemize}

\item[(a)] $T$ is continuous.

\item[(b)] $S(M)\subset T(M).$

\item[(c)] $(S,T)$ is a commuting pair.
\end{itemize}
 Then, $S$ and $T$ have a unique common fixed point.
\end{teo}

\end{document}